\documentclass[12pt, reqno]{amsart}
\usepackage{amsmath,amssymb,amsthm}
\usepackage[english]{babel}
\usepackage[autostyle]{csquotes}
 2

\usepackage[colorlinks=true,linkcolor=blue,citecolor=blue,pdfpagelabels=false]{hyperref}
\newtheorem{thm}{Theorem}[section]
\newtheorem{lem}[thm]{Lemma}
\newtheorem{prop}[thm]{Proposition}

\newtheorem{defn}[thm]{Definition}
\newcommand{\thmref}[1]{Theorem~\ref{#1}}
\newcommand{\lemref}[1]{Lemma~\ref{#1}}

\theoremstyle{remark}

\renewcommand{\geq}{\geqslant}
\renewcommand{\leq}{\leqslant}

\begin{document}

\baselineskip=17pt

\title
    [Simultaneous non-vanishing and sign changes]{Simultaneous non-vanishing and sign changes of Fourier coefficients of
      modular forms}
\author[M. Kumari and R. Murty]{Moni Kumari And M. Ram Murty}
\address{Moni Kumari, School of Mathematical Sciences, National Institute of Science Education and Research, HBNI, Bhubaneswar, Via-Jatni, Khurda, Odisha, 752050, India.}

\email{moni.kumari@niser.ac.in}
  \address{M. Ram Murty, Department of Mathematics and Statistics,
    Queen's University, Kingston, Ontario, K7L 3N6, Canada.}
  \email{murty@queensu.ca}
\date{\today}


\thanks{Research of the second author partially supported by an NSERC Discovery
  Grant.}
\begin{abstract}
  In this article, we give some results on simultaneous non-vanishing and simultaneous sign-changes for the Fourier coefficients of two modular forms. More precisely, given two modular forms $f$ and $g$ with Fourier coefficients $a_n$ and $b_n$ respectively, we consider the following questions:
  existence of infinitely many primes $p$ such that $a_p b_p\neq 0$; simultaneous non-vanishing in the short intervals and in arithmetic progressions;
simultaneous sign changes in short intervals.\\

\end{abstract}

\maketitle

\section{Introduction}
The vanishing or non-vanishing of an arithmetically defined analytic function
is a recurring motif in mathematics.  In recent times, such questions
have arisen in the context of modular forms both of integral weight and
half-integral weight.  
In this paper, we will study simultaneous
non-vanishing of Fourier coefficients of distinct modular forms of integral
weight.  
\par
Throughout we let $k, N$ be positive integers and $p$ be a prime. We write $S_{k}(\Gamma_0(N))$ for the space of cusp forms of weight $k$ for the group $\Gamma_{0}(N).$
Let $f(z)=\sum_{n=1}^{\infty}a_n q^n \in S_{k}(\Gamma_0(N))$ and $g(z)=\sum_{n=1}^{\infty}b_n q^n \in S_{k}(\Gamma_0(N))$ be two non-zero cusp forms which are not a linear combination of CM forms. One of the goals of this paper is to study simultaneous non-vanishing of $a_n, b_n$ partially inspired by a long-standing conjecture of Lehmer which predicts that $\tau(n)\neq 0,~ \mbox{for all}~ n>0,$ where 
$$\Delta(z)=\sum_{n=1}^{\infty} \tau(n)q^n=q \prod_{n=1}^{\infty}(1-q^n)^{24}, ~q:=e^{2 i\pi z}$$
is the unique normalized cusp form of weight $12$ on $SL_2(\mathbb{Z}).$ In relation to Lehmer's conjecture, Serre in his paper \cite{ser}  
motivated the general study of estimating the size of possible gaps in the Fourier expansion of modular forms via the gap function
$$i_f(n):=\mbox{min}\{j\geq 0: a_{n+j}\neq 0\}.$$
He proved that $i_f(n)\ll_f n,$ where $f(z)$ is a cusp form which is not a linear combination of CM forms. In the same paper \cite{ser}, he posed the question of whether one can prove an estimate of the form\\
$$i_f(n)\ll_f n^\delta,$$ where $\delta <1$.
In his paper \cite{kumar}, Kumar Murty first pointed out that
$i_f(n) \ll n^{3/5}$ follows immediately from the celebrated work
of Rankin \cite{rankin} and Selberg \cite{selberg} done in 1939/40.
After that many authors improved the value of $\delta$ (for more detail see \cite{das}).
\par
In the case of level 1
and $f$ an eigenform, Das and Ganguly \cite{das} discovered a clever
argument to show $i_f(n)\ll n^{1/4}$ by combining a classical result of
Bambah and Chowla \cite{bambah} with a congruence of Hatada \cite{hatada}
along with a basic lemma of Murty and Murty \cite{murty-murty}.  
Here is a synopsis of their elegant proof.
In 1947, Bambah and Chowla showed using an elementary argument
that in any interval of length $x^{1/4}$ there is a number $n$ (say)
which can be written
as a sum of two squares.  As $f$ is an eigenform, $a_n$ is
multiplicative.  Hatada's theorem \cite{hatada} implies that
$a_p\equiv 2$ (mod 4), for $p\equiv 1$ (mod 4) and $a_{p^r}\equiv
1 $ (mod 4) if $r$ is even and $p\equiv 3$ (mod 4).  The lemma in
\cite{murty-murty} shows that $a_{p^r} \neq 0$ for $p\equiv 1$ (mod 4)
provided $p$ is sufficiently large.
These congruences
combined with the classical theorem
about factorization of natural numbers that
can be written as a sum of two integral squares now imply $a_n\neq 0$
provided $n$ is coprime to a given finite set of primes.
Thus, one now needs the Bambah-Chowla
theorem with $n$ coprime to a finite set of primes.  One can tweak
the argument in \cite{bambah} to accomodate this extra condition and thus
deduce the non-vanishing result as done in \cite{das}.
Actually, the argument of Bambah and Chowla can be generalized
with considerable latitude. We prove the following
which is of independent interest.
\begin{thm}\label{0}
  Let $r$ and $s$ be natural numbers and set $\alpha = (r-1)(s-1)/rs$.
  There is an effectively computable $C$ (depending only
  on $r$ and $s$) such that in any interval of the form
  $[n, n+ Cn^\alpha]$, there is a number $m$ which can be
  written as
  $$m = A^r + B^s , $$
  with $A$ and $B$ integers.

\end{thm}  
\par
We hasten to highlight that the argument of Das and Ganguly
allows for simultaneous
non-vanishing.  In fact, if $f_1, ..., f_r$ are normalized eigenforms
of level 1, with corresponding Fourier coefficients $a_n(f_j)$,
then one can find an $i$ with $i\ll n^{1/4}$ such that
$$a_{n+i}(f_j) \neq 0, \qquad \forall \quad 1\leq j\leq r. $$
\par
It has been suggested that perhaps $i_f(n) \ll n^{\epsilon}$ for
any $\epsilon >0$.  Perhaps even the stronger conjecture
$i_f(n) \ll 1$ is true (see for example, \cite{kumar}).
\par
It would be nice to extend these results to higher levels
but as the authors in \cite{das} remark, one needs to extend Hatada's
result, which probably can be done, but would take us in a direction
orthogonal to the methods of this paper.  We expect to return to this
question at a later time.
\par
In this paper, we  introduce the analogous concept of gap function $i_{f,g}$ for simultaneous non-vanishing and then we derive a bound for $i_{f,g}$ as small as possible, based on
current knowledge.  One can, of course, consider more general gap functions
for several modular forms.

\begin{thm}\label{1}
Let $f(z)=\sum_{n=1}^{\infty} a_n n^{\frac{k-1}{2}}q^n \in S_{k}(\Gamma_0(N))$ and $g(z)=\sum_{n=1}^{\infty} b_n n^{\frac{k-1}{2}}q^n\in S_{k}(\Gamma_0(N))$ be two newforms which are not CM forms, then there exist infinitely many primes $p$ such that 
$$a_p b_p\neq 0.$$
\end{thm}
\par
Actually, as we show below, the theorem is true without the constraint
that the forms are not CM.  It should be remembered
that the recent solution of the Sato-Tate conjecture for two
distinct eigenforms (see \cite{murty-pujahari} and the references
therein), the theorem is immediate.  This is because there is a positive
density of primes $p$ such that both $a_p$ and $b_p$ are simultaneously
non-zero since the joint Sato-Tate distribution holds for two eigenforms
$f$ and $g$.  
But this is invoking a ``sledgehammer''
result and we underscore that our methods do not make use of this
major advance.  This comment amplifies that there are other  
ways of approaching such questions.  
\par
Now, for $n \in \mathbb{N}$ 
define $$i_{f,g}(n):= \mbox{min}\{m\geq 0: a_{n+m}b_{n+m}\neq 0\}, $$ which is well-defined from the above theorem. We are interested to find the growth of the function $i_{f,g}(n)$ as $n\rightarrow \infty.$ 
In 2014, Lu \cite{lu} by using the result of Chandrasekharan and Narasimhan \cite{cha} proved the following.
$$\sum_{n\leq x}a_n^2 b_n^2= cx +O(x^{\frac{7}{8}+\varepsilon}),$$
where $c$ is a non-zero constant.
It then follows that 
$$i_{f,g}(n)\ll n^{\frac{7}{8}+\varepsilon}.$$
In the present paper we give a better estimate than above.
\begin{thm}\label{2}
Suppose that $f(z)=\sum_{n=1}^{\infty} a_n n^{\frac{k-1}{2}}q^n \in S_{k}(\Gamma_0(N))$ and $g(z)=\sum_{n=1}^{\infty} b_n n^{\frac{k-1}{2}}q^n \in S_{k}(\Gamma_0(N))$ are two newforms with $k>2$ which are not a linear combination of CM forms. Then the following results hold.
\begin{itemize}
\item[(i)]
 For every $\varepsilon>0$, $x>x_0(f,g,\varepsilon)$ and $x^{\frac{7}{17}+\varepsilon}\leq y$ we have 
 \begin{equation}
   |\{ x<n<x+y: a_n b_n\neq 0 \}| \gg_{f,g,\varepsilon} y.
 \end{equation}
In particular, we get that $i_{f,g}(n)\ll_{f,g,\varepsilon} n^{\frac{7}{17}+\varepsilon}.$
\item[(ii)]
For every $\varepsilon>0,~ x\geqslant x_{0}(f,g,\varepsilon),~ y\geqslant x^{\frac{17}{38}+100 \varepsilon}$ and $1\leq a\leq q\leq
x^{\varepsilon}$ with $(a,q)=1,$ we have
\begin{equation}
 |\{ x< n\leq x+y: ~n \equiv a\pmod q~ \mbox{and}~ a_n b_n\neq 0\}|\gg_{f,g,\varepsilon} y/q.
\end{equation} 
\end{itemize}
\end{thm}
\bigskip
In 2009, Kohnen and Sengupta \cite{koh} considered a problem related with the simultaneous sign changes. They proved that,
given two normalized cusp forms $f$ and $g$ of the same level and different weights with totally real algebraic Fourier coefficients, there exists a Galois automorphism $ \sigma$ such that $f^{\sigma}$ and $g^{\sigma}$ have infinitely many Fourier coefficients of the opposite sign. Recently Gun, Kohnen and Rath \cite{gun} removed the dependency on the Galois conjugacy. In fact, they extended their result to arbitrary cusp forms with arbitrary real Fourier coefficients but they assumed that both $f$ and $g$ should have first Fourier coefficient to be non-zero.
 More precisely, they proved the following.
\begin{thm}\label{Sgun}
Let 
\begin{equation*}
f(z)= \sum_{n=1}^{\infty}a_n q^n ~\mbox{and}~  g(z)= \sum_{n=1}^{\infty}b_n q^n
\end{equation*}
be non-zero cusp forms of level $N$ and weights $1<k_1<k_2$ respectively. Suppose that $a_n,b_n$ are real numbers. If $a_1b_1\neq0,$ then there exist infinitely many $n$ such that $a_n b_n>0$ and infinitely many $n$ such that $a_n b_n<0.$
\end{thm}
\bigskip
In this paper we extend the above result by removing the assumption $a_1 b_1\neq 0.$
We prove the following result.
\begin{thm}\label{3}
Let 
\begin{equation*}
f(z)= \sum_{n\geq 1}a_n q^n ~\mbox{and}~  g(z)= \sum_{n\geq 1}b_n q^n
\end{equation*}
be non-zero cusp forms of level $N$ and weights $1<k_1<k_2$ respectively. Further, let $a_n,b_n$ be real numbers. Then there exist infinitely many $n$ such that $a_n b_n>0$ and infinitely many $n$ such that $a_n b_n<0.$
\end{thm}
If $f$ and $g$ are newforms then we have the following quantitative result for the simultaneous sign changes.
\begin{thm}\label{4}
Let $k\geq 2$ be an integer.  Assume that 
\begin{equation*}
f(z)= \sum_{n\geq 1}a_n n^{\frac{k-1}{2}}q^n ~\mbox{and}~  g(z)= \sum_{n\geq 1}b_nn^{\frac{k-1}{2}} q^n
\end{equation*}
are two distinct newforms of weight $k$ on $\Gamma_0({N}).$ Further, let $a_n, b_n$ be real numbers, then for any  $\delta> \frac{7}{8},$ the sequence 
$\{a_n b_n\}_{n\in \mathbb{N}}$ has at least one sign change for $n\in (x,x+x^\delta]$ for sufficiently large $x.$ In particular, the number of sign changes for $n\leq x$ is $\gg x^{1-\delta}.$
\end{thm}
\section{Preliminaries}

In this section, we collect various results from the literature that will be needed in our proofs.
In 1982, Serre \cite[p.174, Cor.2]{ser} in his very famous paper, proved the following result.
\begin{lem}\label{ser}
 Let $f(z)=\sum_{n=1}^{\infty}a_n q^n \in S_k(\Gamma_0(N))$ be a newform with weight $k\geq 2$ which does not have complex multiplication. For every $\epsilon >0$ we have
$$ |\{ p \leqslant x: a_p=0\}|\ll_{f,\epsilon} \frac{x}{(\log x)^{\frac{3}{2}-\epsilon}}.$$
\end{lem}
\bigskip

To prove \thmref{2}, we shall use the concept of $\mathcal{B}$-free numbers which was introduced by Erd\"os in 1966 and later many authors studied the distribution of $\mathcal{B}$-free numbers.
\begin{defn}
($\mathcal{B}$-free numbers): Let $\mathcal{B}=\{b_i: 1<b_1<b_2<...\}$ be a sequence of mutually coprime positive integers for which $\sum_{i=1}^{\infty}\frac{1}{b_i}<\infty.$ A positive integer $n$ is called $ \mathcal{B}$-free if it is not divisible by any element in $\mathcal{B}.$
\end{defn}
By using sieve theory and estimates for multiple exponential sums, Chen and Wu \cite{wu}, studied the distribution of $\mathcal{B}$-free numbers in short intervals as well as in an arithmetic progression and they proved the following results.
\begin{prop} Let $\mathcal{B}$ be a sequence of positive integers satisfying the conditions in the definition of $\mathcal{B}$-free numbers.  Then,\\
\begin{itemize}
\item[(i)]
 for any $\varepsilon > 0,~ x>x_o(\mathcal{B},\varepsilon)$ and $y\geqslant x^{\frac{7}{17}+\varepsilon},$ we have
\begin{equation}\label{wu1}
|\{x<n \leq x+y:n ~\mbox{is} ~\mathcal{B}~\mbox{-free}~\}|\gg_{\mathcal{B},\varepsilon}y,
\end{equation}

\item[(ii)]
 for any $\varepsilon > 0,~ x>x_o(\mathcal{B},\varepsilon)$ and $y\geqslant x^{\frac{17}{38}+100 \varepsilon},~ 1\leq a \leq q\leq x^\varepsilon$ with $((a,q),b)=1,~ \mbox{for all}~ b\in \mathcal{B},$ we have
\begin{equation}\label{wu2}
|\{x<n \leq x+y:n\equiv a \pmod q~\mbox{and}~ n ~\mbox{is} ~\mathcal{B}~\mbox{-free}~\}|\gg_{\mathcal{B},\varepsilon}y/q.
\end{equation}
\end{itemize}
Here the implied constants depend only on $\mathcal{B}$ and $\varepsilon.$
\end{prop}
\bigskip
In the proof of \thmref{3} we use the following theorem of Pribitkin \cite{pri}.
\begin{thm}\label{pri}
Let $F(s)=\sum_{n=1}^{\infty}a_{n}e^{-s\lambda_{n}}$ be a non-trivial general Dirichlet series which converges somewhere, where the sequence $\{a_n\}_{n=1}^{\infty}$ is complex, and the exponent sequence $\{\lambda_n\}_{n=1}^{\infty}$ is real and strictly 
increasing to $ \infty.$ If the function $F$ is holomorphic on the whole real line and has infinitely many real zeros, then there exist infinitely many $n \in \mathbb{N}$ such that $a_{n}>0$ and there exist infinitely many $n\in \mathbb{N}$ such that $a_{n}<0.$
\end{thm}

\section{Proof of \thmref{0}}
We essentially follow Bambah and Chowla \cite{bambah}
and modify their argument to our setting.
Let $t= [n^{1/s}]= n^{1/s} - \theta $ with $0 \leq \theta < 1$.
Let $x_1, x_2$ be positive real numbers such that
$$x_1^r + t^s = n, $$
$$x_2^r + t^s = n + Cn^\alpha $$
with $C$ to be chosen later.
Thus,
$x_2^r - x_1^r = Cn^\alpha. $
Now,
$$x_1 = (n-t^s)^{1/r} \ll n^{(s-1)/rs}, \qquad x_2 \ll n^{(s-1)/rs}, $$
by a simple application of the binomial theorem.
Hence,
$$x_2^{r-1} +x_2^{r-2}x_1+\cdots + x_1^{r-1} \ll n^{(s-1)(r-1)/rs} =  n^\alpha. $$
Now writing
$$(x_2 - x_1) (x_2^{r-1} + \cdots + x_1^{r-1})=x_2^r - x_1^r = Cn^{\alpha}, $$
we immediately see that
$$x_2 - x_1 > 1, $$
for a suitable choice of $C$. (In fact, $C= 2^{rs}rs$ will work.)
Therefore, there is a natural number $N$  in the interval $[x_1, x_2]$
so that
$$n = x_1^r + t^s < N^r + t^s < x_2^r + t^s = n + Cn^\alpha, $$
as desired.
This completes the proof of Theorem \ref{0}.
\par
We remark that there are several variations of this theorem
that can be derived from this proof.  For example, if $f(x)$ is
a monotonic, continuous function for $x$ sufficiently large,
and $f(x) \asymp x^r $, then there is a natural number $m$
such that $m = f(A) + B^s$ for some natural numbers $A,B$ and
with $m\in [n, n+Cn^{\alpha}]$.  In particular, this can be applied
to the norm form $a^2 + Db^2$, with $D$ squarefree.
We record these remarks with the view that the result may have potential applications in other contexts.

\section{Proof of \thmref{1}}
From \lemref{ser}, we have
$$ |\{ p \leqslant x : a_p=0\}|\ll_{f,\epsilon} \frac{x}{(\log x)^{\frac{3}{2}-\epsilon}},$$
and $$ |\{ p \leqslant x : b_p=0\}|\ll_{g,\epsilon} \frac{x}{(\log x)^{\frac{3}{2}-\epsilon}}.$$

Since $a_p b_p=0,$ we have either $a_p=0 ~~\mbox{or}~~ b_p=0.$ 
Hence $$  |\{ p\leqslant x : a_p b_p=0\}|\ll_{f,g,\epsilon} \frac{x}{(\log x)^{\frac{3}{2}-\epsilon}}.$$
By the prime number theorem, we have
 $$\pi(x):=|\{p\leq x \}|\sim \frac{x}{\log x}.$$
Hence
 $$|\{ p\leqslant x: a_p b_p\neq0\}|= \pi(x)- |\{ p\leq x: a_p b_p=0\}|\sim  \frac{x}{\log x}.$$
Thus there exist infinitely many primes $p$ such that
$$a_p b_p\neq0.$$
We make some remarks in the case that either $f$ or $g$ is of CM type.  Suppose first that
$f$ has CM by an order in an imaginary quadratic field $K$ and $g$ does not.
Then, for primes $p$ coprime to the level of $f$, $a_p=0$ if and only if $p$ is inert in $K$.
The density of such primes is $1/2$ and so
$$|\{ p\leqslant x: a_p b_p\neq0\}|= \pi(x)- |\{ p\leq x: a_p b_p=0\}|\gtrsim \frac{x}{2\log x}.$$
Hence, in this case also, there are infinitely many primes $p$ such that $a_pb_p \neq 0$.
If both $f$ and $g$ have CM by two imaginary quadratic fields $K_1, K_2$ (say, respectively),
then we need only choose primes $p$ coprime to the level which split in $K_1$
and $K_2$.  This density is either $1/2$ (if $K_1=K_2$) or $1/4$ (if $K_1\neq K_2$).
Thus, in all cases, \thmref{1} is valid in general.
\section{Proof of \thmref{2}}
Let $S=\{ p: a_p b_p =0\} \cup \{p|N\}.$ Put 
$ \mathcal{B}= S \cup \{p^2: p \notin S \}.$
Clearly $\mathcal{B}$ is a sequence of mutually coprime integers and if $n$ is $\mathcal{B}$-free, then $n$ is square-free and 
$a_n b_n\neq 0$ by using the multiplicative properties of $a_n$ and $b_n.$ Thus \eqref{wu1} and \eqref{wu2} imply the first and second assertions of Theorem 1.2 respectively, if we can show that
 $\sum_{p\in \mathcal{B}}\frac{1}{p}< \infty.$
Since $\sum_{p} 1/p^2 < \infty$, 
it suffices to show that  $$\sum_{p\in S}\frac{1}{p}< \infty.$$
We know, from \lemref{ser} that
\begin{equation*}
\sum_{\substack{p\leqslant x\\ p\in S}}1\ll_{f,g} \frac{x}{(\log x)^{1+\delta}},~~ \mbox{for some}~~ \delta>0.
\end{equation*}
Hence, by partial summation formula, we have 
\begin{align*}
& \sum _{\substack{p\leqslant x\\ p\in S}}\frac{1}{p}
=\frac{1}{x}\sum_{\substack{p\leqslant x\\ p\in S}}1 + \int_{2}^{x} \frac{1}{t^2}(\sum_{\substack{p\leqslant t\\ p\in S}}1) dt
\ll _{f,g} \frac{1}{(\log x)^{1+\delta}}+ \int_{2}^{x}\frac{dt}{t (\log t)^{1+\delta}}\\
&\ll _{f,g} 1.
\end{align*}
This completes the proof of \thmref{2}.

\section{Proof of \thmref{3}}
We assume that either $a_1=0$ or $b_1=0,$ since otherwise by using \thmref{Sgun}, we get the result.
We will show that there exists infinitely many $n\in \mathbb{N}$ such that $a_{n}b_{n}<0$
the other case being similar.
Suppose not, then there exist $n_{0}\in \mathbb{N}$ such that 
$$ a_{n}b_{n}\geq 0,$$ for all $n>n_{0}.$
Set $M=\prod_{p\leq n_{0}}p.$
Clearly, $a_{n}b_{n}\geq0 $ whenever $(n,M)=1$ by our assumption. 
Let
$$ f_1(z):=\sum_{\substack{n\geq 1 \\ (n,M)=1}}a_{n}q^n~~~\mbox{and}~~g_1(z):=\sum_{\substack{n\geq1\\ (n,M)=1}}b_{n}q^n.$$
Then $f_{1}$ and $g_{1}$ are cusp forms of level $NM^2$ and weights $k_{1}$ and $k_{2}$ respectively.
For $s\in \mathbb{C}$ with $\rm Re(s)\gg 1,$ the Rankin-Selberg $L$-function attached to $f_{1}$ and $g_{1}$ is defined by

$$ R_{f_1,g_1}(s):= \sum_{\substack{n\geq1 \\ (n,M)=1}} \frac{a_{n} b_{n}}{n^s}.$$
For  $\rm Re(s)\gg 1,$ set
 \begin{align}\label{Li}
 L_{f_1,g_1}(s)&:= \prod_{p|NM^2}(1-p^{-(2s-(k_{1}+k_{2})+2)}) \zeta(2s-(k_{1}+k_{2})+2)R_{f_1,g_1}(s)\\
               &:=\sum_{n=1}^{\infty}\frac{c_{n}}{n^s}.
              \end{align}

Winnie Li \cite{li} proved that
$ (2\pi)^{-2s}\Gamma(s)\Gamma(s-k_{1}+1) L_{f_1,g_1}(s)$ is entire and we also know that $\Gamma(s)\Gamma(s-k_{1}+1)$ does not have any zeros. Hence $L_{f_1, g_1}(s)$ is an entire function on $\mathbb{C}.$
Let us observe that
the coefficients of this Dirichlet series are non-negative because the
term
$$\prod_{p|NM^2}(1-p^{-(2s-(k_{1}+k_{2})+2)}) \zeta(2s-(k_{1}+k_{2})+2)$$
is the Riemann zeta function with the Euler factors at primes $p|NM^2$ removed
and so is a Dirichlet series with non-negative coefficients.  
Hence $\sum_{n=1}^{\infty}\frac{c_{n}}{n^s}$ is entire with $c_{n}\geqslant 0$ for all $n.$
Now $\sum_{n=1}^{\infty}\frac{c_{n}}{n^s}$ has infinitely many real zeros coming from the real simple poles of the $\Gamma$-function. Then by \thmref{pri},
there exist infinitely many $n$ such that $ c_{n}>0$ and there exist infinitely many $n$ such that $c_{n}<0$ which is a contradiction because $c_{n}\geq 0$ for all $n \in\mathbb{N}.$ This completes our proof.\\
\section{Proof of \thmref{4}}
Recently Meher and the second author \cite{jab}, gave a general criteria for the sign changes of any sequence of real numbers $\{a_n\}_{n\in \mathbb{N}}.$ More precisely, they proved the following.
\begin{thm}\label{murty}
Let $\{a_n\}_{n\in \mathbb{N}}$ be a sequence of real numbers such that
\begin{itemize}
\item[(i)]
$a_n= O(n^\alpha),$
\item[(ii)]
$\sum_{n\leq x}a_n= O(n^\beta),$
\item[(iii)]
$\sum_{n\leq x}a_n^2= cx+O(n^\gamma),$
\end{itemize}
with $\alpha,\beta,\gamma,c\geq 0.$ If $\alpha+\beta < 1,$ then for any $\delta$ satisfying 
$$ \rm{max}\{\alpha+\beta,\gamma\}<\delta<1,$$
the sequence $\{a_n\}_{n \in \mathbb{N}}$ has at least one sign change for $n\in [x,x+x^\delta].$ Consequently, the number of sign changes of $a_n$ for $n\leq x$ is $\gg x^{1-\delta}$ for sufficiently large $x.$
\end{thm}
We will prove \thmref{4}, as an application of the above theorem, for which we have to analyse the stated conditions for the sequence $\{a_n b_n\}_{n\in \mathbb{N}}.$
\begin{itemize}
\item[(i)] Ramanujan-Deligne: 
\begin{equation}
a_n b_n= O( n^\varepsilon) ~\mbox{ for all}~ \varepsilon>0.
\end{equation}
From the paper of Lu \cite{lu}, one can deduce the following results
\item[(ii)]
\begin{equation}
\sum_{n\leq x}a_n b_n\ll x^\frac{3}{5} (\log x)^{-\frac{2\theta}{3}}.\\
\end{equation}
where $\theta= 0.1512....$
\item[(iii)]
$$\sum_{n\leq x}a_n^2 b_n^2= cx +O(x^{\frac{7}{8}+\varepsilon}).$$
\end{itemize}
Hence from \thmref{murty}, we immediately deduce Theorem \ref{4}.

\section{Concluding remarks}

As mentioned earlier, it would be interesting to extend Hatada's congruence
to modular forms of higher level.  This is a research problem of
independent interest and is accessible since there have been
significant advances in the theory of congruences of modular forms.  If one
assumes standard conjectures about distribution of primes such as
Cram\'er's conjecture, then it is easy to deduce that $i_f(n) = O(\log^2 n)$.
The other problem that suggests itself is to obtain estimates with
their dependence on level and weight made explicit.  An initiation
into such an enterprise can be found in the methods of \cite{murty-1} and
\cite{murty-2}.
\par
The analogues of these questions for modular forms of half-integral
weight takes us into a parallel universe of ideas.  There is, of course,
a link between these two worlds provided by Waldspurger's theorem
and the question is equivalent to the simultaneous non-vanishing
of quadratic twists of $L$-series attached to modular forms.  A modest
beginning in this line of research was initiated in \cite{lakshmi}.
$$\quad $$
\noindent{\bf Acknowledgements.}  We thank Soumya Das, Satadal Ganguly,
Jaban Meher and Brundaban Sahu for their comments on an earlier
version of this paper. We also thank the referee for helpful remarks and  suggestions.

\end{document}